# A Family of Median Based Estimators in Simple Random Sampling

[1]Hemant K.Verma, †[1]Rajesh Singh and [2]Florentin Smarandache
Department of Statistics, Banaras Hindu University
Varanasi-221005, India
[2]Chair of Department of Mathematics, University of New Mexico, Gallup, USA

†Corresponding author

**Abstract**

In this paper we have proposed a median based estimator using known value of some population parameter(s) in simple random sampling. Various existing estimators are shown particular members of the proposed estimator. The bias and mean squared error of the proposed estimator is obtained up to the first order of approximation under simple random sampling without replacement. An empirical study is carried out to judge the superiority of proposed estimator over others.

## 1. Introduction

Consider a finite population $U = \{U_1, U_2, ..., U_N\}$ of N distinct and identifiable units. Let Y be the study variable with value $Y_i$ measured on $U_i, i = 1,2,3..., N$. The problem is to estimate the population mean $\overline{Y} = \frac{1}{N}\sum_{i=1}^{N} Y_i$. The simplest estimator of a finite population mean is the sample mean obtained from the simple random sampling without replacement, when there is no auxiliary information available. Sometimes there exists an auxiliary variable X which is positively correlated with the study variable Y. The information available on the auxiliary variable X may be utilized to obtain an efficient estimator of the population mean. The sampling theory describes a wide variety of techniques for using auxiliary information to obtain more efficient estimators. The ratio estimator and the regression estimator are the two important estimators available in the literature which are using the auxiliary information. To know more about the ratio and regression estimators and other related results one may refer to [1-13].

When the population parameters of the auxiliary variable X such as population mean, coefficient of variation, kurtosis, skewness and median are known, a number of modified ratio estimators are proposed in the literature, by extending the usual ratio and Exponential- ratio type estimators.

Before discussing further about the modified ratio estimators and the proposed median based modified ratio estimators the notations and formulae to be used in this paper are described below:

- N - Population size
- n - Sample size
- Y - Study variable
- X - Auxiliary variable
- $\beta_1 = \frac{\mu_3^2}{\mu_2^3}$ Where $\mu_r = \frac{1}{N}\sum_{i=1}^{N}(x_i - \overline{X})^r$, Coefficient of skewness of the auxiliary variable
- $\rho$ - Correlation Co-efficient between X and Y
- $\overline{X}, \overline{Y}$ - Population means
- $\overline{x}, \overline{y}$ - Sample means
- $\overline{M}$, - Average of sample medians of Y
- m - Sample median of Y
- $\beta$ - Regression coefficient of Y on X
- B (.) - Bias of the estimator
- V (.) - Variance of the estimator
- MSE (.) - Mean squared error of the estimator
- $PRE(e,p) = \frac{MSE(e)}{MSE(e)} * 100$ - Percentage relative efficiency of the proposed estimator p with respect to the existing estimator e.

The formulae for computing various measures including the variance and the covariance of the SRSWOR sample mean and sample median are as follows:

$$V(\overline{y}) = \frac{1}{^{N}C_n}\sum_{i=1}^{^{N}C_n}(y_i - \overline{Y})^2 = \frac{1-f}{n}S_y^2, V(\overline{x}) = \frac{1}{^{N}C_n}\sum_{i=1}^{^{N}C_n}(x_i - \overline{X})^2 = \frac{1-f}{n}S_x^2, V(m) = \frac{1}{^{N}C_n}\sum_{i=1}^{^{N}C_n}(m_i - \overline{M})^2 \ ,$$

$$Cov(\overline{y},\overline{x}) = \frac{1}{^{N}C_n}\sum_{i=1}^{^{N}C_n}(x_i - \overline{X})(y_i - \overline{Y}) = \frac{1-f}{n}\frac{1}{N-1}\sum_{i=1}^{N}(x_i - \overline{X})(y_i - \overline{Y}),$$

$$Cov(\overline{y},m) = \frac{1}{^{N}C_n}\sum_{i=1}^{^{N}C_n}(m_i - \overline{M})(y_i - \overline{Y}),$$

$$C'_{xx} = \frac{V(\overline{x})}{\overline{X}^2}, C'_{mm} = \frac{V(m)}{\overline{M}^2}, C'_{ym} = \frac{Cov(\overline{y},m)}{\overline{M}\overline{Y}}, C'_{yx} = \frac{Cov(\overline{y},\overline{x})}{\overline{X}\overline{Y}}$$

Where $f = \frac{n}{N}; S_y^2 = \frac{1}{N-1}\sum_{i=1}^{N}(y_i - \bar{Y})^2, S_x^2 = \frac{1}{N-1}\sum_{i=1}^{N}(x_i - \bar{X})^2$ ,

In the case of simple random sampling without replacement (SRSWOR), the sample mean $\bar{y}$ is used to estimate the population mean $\bar{Y}$ . That is the estimator of $\bar{Y} = \bar{Y}_r = \bar{y}$ with the variance

$$V(\bar{Y}_r) = \frac{1-f}{n} S_y^2 \quad (1.1)$$

The classical Ratio estimator for estimating the population mean $\bar{Y}$ of the study variable Y is defined as $\bar{Y}_R = \frac{\bar{y}}{\bar{x}}\bar{X}$ . The bias and mean squared error of $\bar{Y}_R$ are given as below:

$$B(\bar{Y}_R) = \bar{Y}\{C'_{xx} - C'_{yx}\} \quad (1.2)$$

$$MSE(\bar{Y}_R) = V(\bar{y}) + R'^2 V(\bar{x}) - 2R' Cov(\bar{y}, \bar{x}) \qquad \text{where } R' = \frac{\bar{Y}}{\bar{X}} \quad (1.3)$$

**2. Proposed estimator**

Suppose

$$t_0 = \bar{y}, \quad t_1 = \bar{y}\left[\frac{\bar{M}^*}{\alpha m^* + (1-\alpha)\bar{M}^*}\right]^g, \quad t_2 = \bar{y}\exp\left[\frac{\delta(\bar{M}^* - m^*)}{\bar{M}^* + m^*}\right] \quad \text{where } \bar{M}^* = a\bar{M} + b,\ m^* = am + b$$

Such that $t_0,\ t_1,\ t_2 \in w$, where $w$ denotes the set of all possible ratio type estimators for estimating the population mean $\bar{Y}$. By definition the set $w$ is a linear variety, if

$$t = w_0\bar{y} + w_1 t_1 + w_2 t_2 \qquad \in W, \quad (2.1)$$

$$\text{for} \sum_{i=0}^{2} w_i = 1 \qquad w_i \in R \quad (2.2)$$

where $w_i$ (i=0, 1, 2) denotes the statistical constants and R denotes the set of real numbers.

$$\text{Also, } t_1 = \bar{y}\left[\frac{\bar{M}^*}{\alpha m^* + (1-\alpha)\bar{M}^*}\right]^g, \quad t_2 = \bar{y}\exp\left[\frac{\delta(\bar{M}^* - m^*)}{\bar{M}^* + m^*}\right]$$

and $\bar{M}^* = a\bar{M} + b,\ m^* = am + b.$

To obtain the bias and MSE expressions of the estimator t, we write

$$\bar{y} = \bar{Y}(1 + e_0), \qquad m = \bar{M}(1 + e_1)$$

such that

$$E(e_0) = E(e_1) = 0,$$

$$E(e_0^2)=\frac{V(\bar{y})}{\bar{Y}^2},\quad E(e_1^2)=\frac{V(m)}{\bar{M}^2}=C_{mm},\quad E(e_0e_1)=\frac{Cov(\bar{y},m)}{\bar{Y}\bar{M}}=C_{ym}$$

Expressing the estimator t in terms of e's, we have

$$t=\bar{Y}(1+e_0)\left[w_0+w_1(1+\upsilon\alpha e_1)^{-g}+w_2\exp\left\{\left(-\frac{\upsilon\delta e_1}{2}\right)\left(1+\frac{\upsilon e_1}{2}\right)^{-1}\right\}\right] \tag{2.3}$$

where $\upsilon=\frac{a\bar{M}}{a\bar{M}+b}$.

Expanding the right hand side of equation(2.3) up to the first order of approximation, we get

$$t\cong\bar{Y}\left[1-\upsilon we_1+e_0+\upsilon^2\left(w_1\frac{g(g+1)}{2}\alpha^2+\left(\frac{\delta}{4}-\frac{\delta^2}{8}\right)w_2\right)e_1^2-\upsilon we_0e_1\right] \tag{2.4}$$

where $\quad w=\alpha gw_1+\frac{\delta}{2}w_2.$ (2.5)

Taking expectations of both sides of (2.4) and then subtracting $\bar{Y}$ from both sides, we get the biases of the estimators, up to the first order of approximation as

$$B(t)=\bar{Y}\left[\upsilon^2\left\{w_1\frac{g(g+1)}{2}\alpha^2+\left(\frac{\delta}{4}-\frac{\delta^2}{8}\right)w_2\right\}C_{mm}-\upsilon wC_{ym}\right] \tag{2.6}$$

$$B(t_1)=\bar{Y}g\alpha\upsilon\left[\frac{\alpha\upsilon(g+1)}{2}C_{mm}-C_{ym}\right] \tag{2.7}$$

$$B(t_2)=\bar{Y}\left[\left(\frac{\delta\upsilon^2}{4}+\frac{\delta^2\upsilon^2}{8}\right)C_{mm}-\frac{\delta\upsilon}{2}C_{ym}\right] \tag{2.8}$$

From (2.4), we have

$$t-\bar{Y}\cong\bar{Y}(e_0+\upsilon we_1) \tag{2.9}$$

Squaring both sides of (2.9) and then taking expectations, we get the MSE of the estimator t, up to the first order of approximation as

$$MSE(t)=V(\bar{y})+\upsilon^2R^2w^2V(m)-2\upsilon RwCov(\bar{y},m) \tag{2.10}$$

where $R=\frac{\bar{Y}}{\bar{M}}$.

MSE(t) will be minimum, when

$$w=\frac{1}{\upsilon R}\frac{Cov(\bar{y},m)}{V(m)}=k(say) \tag{2.11}$$

Putting the value of w(=k) in (2.10), we get the minimum MSE of the estimator t, as

$$\min . MSE(t) = V(\bar{y})(1-\rho^2) \tag{2.12}$$

The minimum MSE of the estimator t is same as that of traditional linear regression estimator.
From (2.5) and (2.11), we have

$$\alpha g w_1 + \frac{\delta}{2} w_2 = k \tag{2.13}$$

From (2.2) and (2.13), we have only two equations in three unknowns. It is not possible to find the unique values of $w_{i's}$ (i=0, 1, 2). In order to get unique values for $w_{i's}$, we shall impose the linear restriction

$$w_0 B(\bar{y}) + w_1 B(t_1) + w_2 B(t_2) = 0 \tag{2.14}$$

Equations (2.2), (2.11) and (2.14) can be written in matrix form as

$$\begin{bmatrix} 1 & 1 & 1 \\ 0 & \alpha g & \frac{\delta}{2} \\ 0 & B(t_1) & B(t_2) \end{bmatrix} \begin{bmatrix} w_0 \\ w_1 \\ w_2 \end{bmatrix} = \begin{bmatrix} 1 \\ k \\ 0 \end{bmatrix} \tag{2.15}$$

Using (2.15) we get the unique value of $w_{i's}$ (i=0, 1, 2) as

$$\left.\begin{aligned} w_0 &= \frac{\Delta_0}{\Delta_r} \\ w_1 &= \frac{\Delta_1}{\Delta_r} \\ w_2 &= \frac{\Delta_2}{\Delta_r} \end{aligned}\right\} \quad \text{where} \quad \begin{aligned} \Delta_r &= \alpha g B(t_2) - \frac{\delta}{2} B(t_1) \\ \Delta_0 &= B(t_2)(\alpha g - k) + \frac{1}{2} B(t_1)\left(k - \frac{\delta}{2}\right) \\ \Delta_1 &= k B(t_2) \\ \Delta_2 &= -k B(t_1) \end{aligned} \tag{2.16}$$

**Table 2.1: Some members of the proposed estimator**

| $w_0$ | $w_1$ | $w_2$ | **a** | **b** | $\alpha$ | **g** | $\delta$ | **Estimators** |
|---|---|---|---|---|---|---|---|---|
| 1 | 0 | 0 | - | - | - | - | - | $q_1 = \bar{y}$ |
| 0 | 1 | 0 | 1 | 0 | 1 | 1 | - | $q_2 = \bar{y}\frac{\bar{M}}{m}$ |
| 0 | 1 | 0 | $\beta_1$ | $\rho$ | 1 | 1 | - | $q_3 = \bar{y}\left[\frac{\beta_1\bar{M}+\rho}{\beta_1 m+\rho}\right]$ |
| 0 | 1 | 0 | $\rho$ | $\beta_1$ | 1 | 1 | - | $q_4 = \bar{y}\left[\frac{\rho\bar{M}+\beta_1}{\rho m+\beta_1}\right]$ |
| 0 | 0 | 1 | 1 | 0 | - | - | 1 | $q_5 = \bar{y}\exp\left[\frac{(\bar{M}-m)}{\bar{M}+m}\right]$ |
| 0 | 0 | 1 | $\beta_1$ | $\rho$ | - | - | 1 | $q_6 = \bar{y}\exp\left[\frac{\beta_1(\bar{M}-m)}{\beta_1(\bar{M}+m)+2\rho}\right]$ |
| 0 | 0 | 1 | $\rho$ | $\beta_1$ | - | - | 1 | $q_7 = \bar{y}\exp\left[\frac{\rho(\bar{M}-m)}{\rho(\bar{M}+m)+2\beta_1}\right]$ |
| 0 | 1 | 1 | $\beta_1$ | $\rho$ | 1 | 1 | 1 | $q_8 = \bar{y}\left[\frac{\beta_1\bar{M}+\rho}{\beta_1 m+\rho}\right]+\bar{y}\exp\left[\frac{\beta_1(\bar{M}-m)}{\beta_1(\bar{M}+m)+2\rho}\right]$ |
| 0 | 1 | 1 | $\rho$ | $\beta_1$ | 1 | 1 | 1 | $q_9 = \bar{y}\left[\frac{\rho\bar{M}+\beta_1}{\rho m+\beta_1}\right]+\bar{y}\exp\left[\frac{\rho(\bar{M}-m)}{\rho(\bar{M}+m)+2\beta_1}\right]$ |
| 0 | 1 | 1 | 1 | 0 | 1 | 1 | 1 | $q_{10} = \bar{y}\frac{\bar{M}}{m}+\bar{y}\exp\left[\frac{(\bar{M}-m)}{\bar{M}+m}\right]$ |

## 3. Empirical Study

For numerical illustration we consider: the population 1 and 2 taken from [14] pageno.177, the population 3 is taken from [15] page no.104. The parameter values and constants computed for the above populations are given in the Table 3.1. MSE for the proposed and existing estimators computed for the three populations are given in the Table 3.2 whereas the PRE for the proposed and existing estimators computed for the three populations are given in the Table 3.3.

**Table: 3.1 Parameter values and constants for 3 different populations**

| Parameters | For sample size n=3 | | | For sample size n=5 | | |
|---|---|---|---|---|---|---|
| | Popln-1 | Popln-2 | Popln-3 | Popln-1 | Popln-2 | Popln-3 |
| N | 34 | 34 | 20 | 34 | 34 | 20 |
| n | 3 | 3 | 3 | 5 | 5 | 5 |
| ${}^{N}C_{n}$ | 5984 | 5984 | 1140 | 278256 | 278256 | 15504 |
| $\overline{Y}$ | 856.4118 | 856.4118 | 41.5 | 856.4118 | 856.4118 | 41.5 |
| $\overline{M}$ | 747.7223 | 747.7223 | 40.2351 | 736.9811 | 736.9811 | 40.0552 |
| $\overline{X}$ | 208.8824 | 199.4412 | 441.95 | 208.8824 | 199.4412 | 441.95 |
| $\beta_1$ | 0.8732 | 1.2758 | 1.0694 | 0.8732 | 1.2758 | 1.0694 |
| R | 1.1453 | 1.1453 | 1.0314 | 1.1621 | 1.1621 | 1.0361 |
| $V(\bar{y})$ | 163356.4086 | 163356.4086 | 27.1254 | 91690.3713 | 91690.3713 | 14.3605 |
| $V(\bar{x})$ | 6884.4455 | 6857.8555 | 2894.3089 | 3864.1726 | 3849.248 | 1532.2812 |
| V(m) | 101127.6164 | 101127.6164 | 26.0605 | 58464.8803 | 58464.8803 | 10.6370 |
| $Cov(\bar{y}, m)$ | 90236.2939 | 90236.2939 | 21.0918 | 48074.9542 | 48074.9542 | 9.0665 |
| $Cov(\bar{y}, \bar{x})$ | 15061.4011 | 14905.0488 | 182.7425 | 8453.8187 | 8366.0597 | 96.7461 |
| $\rho$ | 0.4491 | 0.4453 | 0.6522 | 0.4491 | 0.4453 | 0.6522 |

**Table: 3.2. Variance / Mean squared error of the existing and proposed estimators**

| Estimators | For sample size n=3 | | | For sample size n=5 | | |
|---|---|---|---|---|---|---|
| | Population-1 | Population-2 | Population-3 | Population-1 | Population-2 | Population-3 |
| $q_1$ | 163356.41 | 163356.41 | 27.13 | 91690.37 | 91690.37 | 14.36 |
| $q_2$ | 89314.58 | 89314.58 | 11.34 | 58908.17 | 58908.17 | 6.99 |
| $q_3$ | 89274.35 | 89287.26 | 11.17 | 58876.02 | 58886.34 | 6.93 |
| $q_4$ | 89163.43 | 89092.75 | 10.92 | 58787.24 | 58730.58 | 6.85 |
| $q_5$ | 93169.40 | 93169.40 | 12.30 | 55561.98 | 55561.98 | 7.82 |
| $q_6$ | 93194.86 | 93186.68 | 12.42 | 55573.42 | 55569.74 | 7.88 |
| $q_7$ | 93265.64 | 93311.19 | 12.62 | 55605.24 | 55625.75 | 7.97 |
| $q_8$ | 113764.16 | 113810.72 | 21.52 | 76860.57 | 76891.47 | 10.66 |
| $q_9$ | 151049.79 | 150701.09 | 22.00 | 101236.37 | 101004.87 | 10.99 |
| $q_{10}$ | 151791.97 | 151791.97 | 24.24 | 101728.97 | 101728.97 | 11.87 |
| t(opt) | 82838.45 | 82838.45 | 10.05 | 52158.93 | 52158.93 | 6.63 |

**Table: 3.3. Percentage Relative Efficiency of estimators with respect to $\bar{y}$**

| Estimators | For sample size n=3 | | | For sample size n=5 | | |
|---|---|---|---|---|---|---|
| | Population-1 | Population-2 | Population-3 | Population-1 | Population-2 | Population-3 |
| $q_1$ | 100 | 100 | 100 | 100 | 100 | 100 |
| $q_2$ | 182.90 | 182.90 | 239.191236 | 155.65 | 155.65 | 205.40 |
| $q_3$ | 182.98 | 182.96 | 242.877047 | 155.73 | 155.71 | 207.12 |
| $q_4$ | 183.21 | 183.36 | 248.504702 | 155.97 | 156.12 | 209.64 |
| $q_5$ | 175.33 | 175.33 | 220.500742 | 165.02 | 165.02 | 183.60 |
| $q_6$ | 175.28 | 175.30 | 218.381298 | 164.99 | 165.00 | 182.30 |
| $q_7$ | 175.15 | 175.07 | 214.915968 | 164.90 | 164.83 | 180.16 |
| $q_8$ | 143.59 | 143.53 | 126.034732 | 119.29 | 119.25 | 134.70 |
| $q_9$ | 108.15 | 108.40 | 123.254986 | 90.57 | 90.78 | 130.57 |
| $q_{10}$ | 107.62 | 107.62 | 111.896010 | 90.13 | 90.13 | 120.97 |
| t(opt) | 197.20 | 197.20 | 269.771157 | 175.79 | 175.79 | 216.51 |

## 4. Conclusion

From empirical study we conclude that the proposed estimator under optimum conditions perform better than other estimators considered in this paper. The relative efficiencies and MSE of various estimators are listed in Table 3.2 and 3.3.